\documentstyle[11pt,amssymb,amstex]{amsart}

\newtheorem{thm}{Theorem}[section]
\newtheorem{lem}[thm]{Lemma}

\def\a{\alpha}
\def\o{\omega}
\def\O{\Omega}
\def\e{\varepsilon}
\begin{document}

\title[Unified martingale and ergodic theorems with continuous time]
{Vector valued unified martingale and ergodic theorems with
continuous parameter}

\maketitle
\begin{center}
\author{ F.A. Shahidi  \footnote{Email: farruh.shahidi@@gmail.com}  }\\\vspace{12pt}
\address{{\em{Department of Mathematics,\\ The Pennsylvania State University,\\University Park, 16802, PA, USA}}}

\end{center}

\begin{abstract} We prove  martingale-ergodic and
ergodic-martingale theorems with continuous parameter for vector
valued Bochner integrable functions. We first provide almost
everywhere convergence of vector valued martingales with
continuous parameter. The norm as well as almost everywhere
convergence of martingale-ergodic and ergodic-martingale averages
are given. We also obtain dominant and maximal inequalities.
Finally, we show that a.e. martingale-ergodic and
ergodic-martingale theorems will coincide under certain
assumptions.

\end{abstract}

\section{Introduction}
An interesting connection in terms of the behavior and convergence
between two fundamental mathematical objects~--- martingales and
ergodic averages has been known since S. Kakutani \cite{Kak}, who
asked for a possible unification of martingale convergence and
ergodic theorems. Several attempts have been done since then (see
\cite{Kach2} for review and references), but none of them was
comprehensive. Quite recently, A.G. Kachurovskii
\cite{Kach2},\cite{Kach1} solved this problem by defining a
martingale-ergodic processes as the composition of martingales and
ergodic averages. For $f\in L_p, p\ge 1,$ if $f_n=E(f|F_n)$ is a
regular martingale, where $E(\cdot|F)$ is a conditional
expectation operator and $A_mf=\frac1m\sum\limits_{i=0}^{m-1}
T^if,$ where $T$ is an $L_1-L_{\infty}$ contraction, then he
proved the following

\begin{thm}\cite{Kach2},\cite{Kach1}.

\begin{enumerate}
\item \begin{itemize}
    \item[(a)] If $f\in L_p,\ p\ge 1,$ then $E(A_mf|F_n)$
    converges in $L_p$ norm as $n,m\to\infty;$
    \item[(b)] If $f\in L_1$ and $sup_n |E(f|F_n)|$ is integrable, then $E(A_mf|F_n)$
    converges   almost everywhere as $n,m\to\infty$.
\end{itemize}

    \item \begin{itemize}
    \item[(c)] If $f\in L_p,\ p\ge 1,$ then $A_mE(f|F_n)$
    converges in $L_p$ norm as $n,m\to\infty;$
    \item[(d)] If $f\in L_1$ and $sup_m |A_mf|$ is integrable, then $A_mE(f|F_n)$
    converges   almost everywhere as $n,m\to\infty$.
\end{itemize}
\end{enumerate}

\end{thm}
While the first part of this theorem is referred as a
martingale-ergodic theorem, second part is known as
ergodic-martingale theorem. In fact, this theorem puts martingale
convergence and ergodic theorems into one superstructure, from
which both martingale convergence and ergodic theorems can be
obtained as degenerate cases.

The continuous parameter analogue of the above theorem was solved
by I.V. Podvigin as follows

\begin{thm}\cite{pod1}
\begin{enumerate}\item\begin{itemize}
    \item[(a)] If $f\in L_p,\ p\ge 1$ then $E(A_tf|F_s)$
    converges in $L_p$ norm as $t,s\to\infty$;
    \item[(b)] If $f\in L_1$ and $sup_s |E(f|F_s)|$ is integrable, then
    $E(A_tf|F_s)$ converges almost everywhere as $t,s\to\infty$.
\end{itemize}
    \item\begin{itemize}
    \item[(c)] If $f\in L_p,\ p\ge 1,$ then $A_tE(f|F_s)$
    converges in $L_p$ norm as $t,s\to\infty$;
    \item[(d)] If $f\in L_1$ and $sup_t |A_tf|$ is integrable, then
    $A_tE(f|F_s)$ converges almost everywhere as $t,s\to\infty$.

    \end{itemize}

\end{enumerate}
Here $F_s-$ an increasing family of $\sigma-$ subalgebras
$A_tf=\frac1t\int\limits_0^tT_{\tau}fd\tau$ and $\{T_t,t\ge 0\}$
is a semigroup of linear $L_1-L_{\infty}$ contractions.
\end{thm}

Note that there many analogues and generalizations of martingale
convergence and ergodic theorems. For example, vector valued
ergodic theorem for $1-$ parameter semigroup of operators was
given by Sh. Hasegawa, R. Sato and Sh. Tsurumi in \cite{Sato1}.
The result was also extended to multiparameter case under suitable
assumptions in \cite{Sato2}. Related problems are also considered
in \cite{Sato4}. This motivates us to provide the above theorem in
other settings. The purpose of this paper is to give the latter
theorem in vector valued settings. Namely, we prove
martingale-ergodic and ergodic-martingale theorems with continuous
parameter for vector valued Bochner integrable functions. As is
done by \cite{Kach1}, \cite{pod1}, we also prove dominant and
maximal inequalities. We also show that the condition of
integrability of supremum is not necessary under the assumption
that conditional expectation operator and ergodic average commute.
This is the vector valued analogue of the result given in
\cite{pod2} for continuous parameter processes. We also note that
the vector valued analogue of Theorem 1.1 has been considered in
\cite{ShGa}.

To our knowledge, we do not seem to have vector valued a.e.
martingale convergence theorem with continuous parameter. Hence in
the next section we prove this convergence. The main result of the
paper is given in section 3. We use the notation and terminology
as used in \cite{pod1}, \cite{ShGa}.

\section{Preliminaries}

In this section we prove a vector valued martingale convergence
theorem with continuous parameter.

 Throughout this paper by $X$ we mean a reflexive Banach space
with the norm $||\cdot||_X$ and by $(\O ,\beta, \mu)$  a finite
measure space. By $L_p(X)=L_p(\O, X), \ 1\le p<\infty$ we denote
the Banach space of $X$ valued measurable functions $f$ on $\O$
with the norm defined as

$$||f||_p=\left(\int_{\O}||f(\o)||_X^pd\mu\right)^{\frac 1p}.$$

We just write $L_p$ when $X=R.$

Let $\{T_t, t\ge 0\}$ be a flow of linear $L_1-L_{\infty}$
contractions acting in $L_1(\O, X)$. That is, for any $t\ge 0,$

$||T_tf||_1\le ||f||$ and $||T_tf||_{\infty}\le ||f||_{\infty}$,
where
$$||f||_1=\int\limits_{\O}||f(\o)||_X d\mu$$ and
$$||f||_{\infty}=inf\{\lambda: ||f(\o)||_X\le\lambda \ a.e\}.$$

A flow of linear operators $\{T_t, t\ge 0\}$ in $L_1(\O,X)$ is
\textit{  strongly continuous semigroup} if
\begin{itemize}
    \item $T_0=id$
    \item $T_{t_1}T_{t_2}=T_{t_1+t_2}$ for all $t_1,t_2>0$
    \item $\lim\limits_{t_1\to t_2}||T_{t_1}f-T_{t_2}f||_1=0$
    for any $f\in L_1(\O, X)$ and $t_2>0.$
\end{itemize}

Henceforth, $\{T_t, t\ge 0\}$ will be a strongly continuous
semigroup of linear  $L_1-L_{\infty}$ contractions unless
otherwise mentioned.

In \cite{Sato1} it is shown that if $f\in L_p(\O, X), p\ge 1,$
then $\frac1t\int\limits_o^tT_{\tau}f(\o)d\tau\in L_p(\O, X).$ In
this settings, we define the ergodic average as follows
$$A_tf(\o)=\frac1t\int\limits_o^tT_{\tau}f(\o)d\tau,\ \ f\in L_1(\O, X),\ \ t>0.$$

The following theorem is an a.e. convergence theorem for the above
ergodic average.

\begin{thm}\cite{Sato1} Let $X$ be a reflexive Banach space and $\{T_t, t\ge 0\}$ be a strongly continuous
semigroup of linear $L_1-L_{\infty}$ contractions on $L_1(\O, X).$
If $1\le p< \infty$ and $f\in L_p(\O, X),$ then the limit
$$\lim\limits_{t\to\infty}\frac1t\int\limits_o^tT_{\tau}f(\o)d\tau$$
exists for almost all $\o\in\O.$

\end{thm}

It is to note that the above theorems were given for slightly
general type of operators $\{T_t\},$ that is, the operators
$\{T_t\}$ should be contractions with respect to $L_1$ norm and
bounded with respect to $L_{\infty}$ norm.

Let $F$ be a $\sigma-$ algebra and $F_1$ be its $\sigma-$
subalgebra.

\begin{thm} \cite{neveu}
\begin{enumerate}
    \item There exists a linear operator $E(\cdot|F):L_1(\O, X)\rightarrow L_1(\O, X)$ such that
$$\int\limits_BE(f|F)d\mu=\int\limits_Bfd\mu$$
for any $f\in L_1(\O, X)$ and $B\subset F_1.$
    \item For every continuous linear functional $g$ and $f\in L_1(\O, X),$
    the function $g(f)$ is integrable and
    $$g(E(f|F))=E(g(f)|F).$$
\end{enumerate}
\end{thm}

 By $E(f|F)$ we denote the conditional expectation of
$f\in L_p.$ Let $F_s,\ s\in R$ be a family of monotonically
increasing (decreasing) sub-$\sigma-$algebras such that
$F_s\uparrow F_{\infty}$ ($F_s\downarrow F_{\infty}$) as
$s\to\infty.$ Unless otherwise stated, we assume that the family
of sub-$\sigma-$ algebras is increasing. We also keep in mind that
the results in this section, which hold for increasing family also
hold for decreasing family of sub-$\sigma-$algebras. A stochastic
process $f_s$ in $L_p(\O, X), \ 1\le p<\infty$ is said to be an
\textit{ordinary (reversed) martingale} if for all $s_1, s_2\in S$
with $s_1<s_2 (s_1>s_2)$ one has $E(f_{s_2}|F_{s_1})=f_{s_1}.$ A
\textit{regular} martingale is given by $f_s=E(f|F_s),$ where
$f\in L_p(\O, X), \ 1\le p<\infty.$ There is a norm convergence
theorem for vector valued martingales with continuous parameter
\cite{Vaxan}. But, we were not able to find any theorem concerning
a.e. convergence for them.  Below we are going to provide this
convergence.

\begin{lem} Let $\{(g_s^i, s\in R), i\in I\}$ be a countable
family of real valued submartingales such that
$$\sup\limits_{s\in R}\int\sup\limits_{i\in I}(g_s^i)^+d\mu< \infty .$$

Then each submartingale converge a.e. to an integrable limit
$g_{\infty}^i,\ i\in I$ and
$$\sup\limits_{i\in I}g_s^i=\sup\limits_{i\in I}g_{\infty}^i$$
as $s\to\infty.$
\end{lem}

\begin{pf} The condition of the lemma implies that $\sup\limits_{s\in R}\int(g_s^i)^+d\mu$
is finite for all $i\in I,$ therefore by Doob's  convergence
theorem for submartingales (see, for example \cite{Oksendal},
Appendix C) the limits
$$g_{\infty}^i=\lim\limits_{s\to\infty}g_s^i$$
exists a.e. Since $g_s^i$ is a submartingale for all $i\in I,$
then $\sup\limits_{i\in I}g_s^i$ is also a submartingale. Due to
the condition of the lemma and Doob's convergence theorem for
submartingales (see \cite{Oksendal})we conclude again that the
limit

$$g_{\infty}=\lim\limits_{s\to\infty}\sup\limits_{i\in I}(g_s^i)$$
exists a.e. This limit clearly dominates each $g_{\infty}^i(i\in
I)$ and thus also their supremum, i.e.
$g_{\infty}\ge\sup\limits_{i\in I}g_{\infty}^i.$

We will show that $\int g_{\infty}d\mu=\int\sup\limits_{i\in
I}(g_{\infty}^i)d\mu$ in order to show that the above inequality
in fact an equality.

Let $(I_p), p\in N$ be a sequence of finite subsets of $I$
increasing to $I$ as $p\to\infty.$ Then the integral
$\int\sup\limits_{i\in I_p}g_s^id\mu$ clearly increases as $p$
increases. Moreover, it also increases with $s (s\in R)$ since
$(\sup\limits_{i\in I_p}g_s^i,\ s\in R)$ is a submartingale for
every $p.$

Note that the expression

$$S=\sup\limits_{p\in N, s\in R}\int\sup\limits_{i\in I_p}g_s^id\mu=
\sup\limits_{s\in R}\int\sup\limits_{i\in I}g_s^id\mu$$
 is dominated by $\sup\limits_{s\in R}\int\sup\limits_{i\in
 I}(g_s^i)^+d\mu$ and hence is finite. Therefore, for every $\e>0$
there exists at least one pair $p_{\e}\in N,\ s_{\e}\in R^+$ such
that
$$\int\sup\limits_{i\in I_p}g_s^id\mu\ge S-\e$$
if $p=p_{\e},\ s=s_{\e}.$ Since the above supremum increases with
$p$ as well as with $s,$ then the above inequality holds for $p\ge
p_{\e},\ s\ge s_{\e}.$ Note that the function
$g_{\infty}-\sup\limits_{i\in I_p}g_{\infty}^i$ is the limit of
positive sequence of functions $(\sup\limits_{i\in I}g_s^i-
\sup\limits_{i\in I_p}g_s^i,\ s\in R)$ so that Fatou's lemma
implies that

$$\int(g_{\infty}-\sup\limits_{i\in I_p}g_{\infty}^i)d\mu\le\liminf\limits_{s\to\infty}\int((\sup\limits_{i\in I}g_s^i-
\sup\limits_{i\in I_p}g_s^i)d\mu\le S-(S-\e)=\e.$$

Therefore, $\int(g_{\infty}-\sup\limits_{i\in
I_p}g_{\infty}^i)d\mu\le\e$ and so $g_{\infty}=\sup\limits_{i\in
I}g_{\infty}^i.$

\end{pf}
\begin{thm} Let $X$ be a separable Banach space which is the dual
of a separable Banach space and $F_s$ be an increasing family of
sub-$\sigma-$algebras. Then for any $f\in L_1(\O, X)$

$$\lim\limits_{s\to\infty}E(f|F_s)= E(f|F_{\infty})$$
a.e. on $\O.$

\end{thm}
Note that every separable reflexive Banach space satisfies the
condition put on $X.$

\begin{pf}

Firstly, note that for any continuous linear functional $g\in X',$
the sequence $g(E(f|F_s))$ is a martingale as (2) of Theorem 2.2
shows that for $s_2>s_1$
$$g(E(f|F_{s_2}))=E(g(f)|F_{s_2})=E(g(f)|F_{s_1})=g(E(f|F_{s_1})).$$
One can also see that for any $g\in X'$

$$g(E(f|F_s))=E(g(f)|F_s)\rightarrow E(g(f)|F_{\infty})=g(E(f|F_{\infty}))$$
outside a set $\O_g$(which actually depends on $g$) of zero
measure as $s\to\infty$ by convergence of read valued martingale
\cite{Oksendal}.

Now assume that the separable Banach space $X$ is the dual of a
(necessarily) separable space $Y$ and let us identify this space
with the subspace of $X',$ the dual of $X.$ Let us denote by $D$ a
dense subset of unit ball in $Y$ which we can choose countable as
$Y$ is separable. Then the equality $\sup\limits_{g\in
D}g(x)=||x||_X$ holds for all $x\in X. $ Indeed, one can see that
$g(x)\le ||g||||x||_X$ implies $||x||_X\ge\frac{g(x)}{||g||},$ and
so $||x||_X\ge \sup\limits_{g\in D}\frac{g(x)}{||g||}.$ Since
there exist $x_0\in X$ and $g_0\in X'$ such that
$g_0(x_0)=||x_0||_X||g_0||,$ then $||x||_X= \sup\limits_{g\in
D}\frac{g(x)}{||g||}.$

Further, take any fixed $a\in X,$ and consider the countable
family of martingales
$$\{(g(E(f|F_s)-a), s\in R), g\in D\}.$$

Since
$$|g(E(f|F_s))|\le ||E(f|F_s)||_X\le E(||f||_X|F_s)$$
for all $g\in D$ by contraction property of the conditional
expectation, then the above family satisfies the condition of
Lemma 2.3 and hence by applying it we get

$$||E(f|F_s)-a||_X\rightarrow ||E(f|F_{\infty})-a||_X$$
a.e. as $s\to\infty,$ for all $a\in X.$ From this it follows that

$$\mu\{\lim\limits_{s\to\infty}||E(f|F_s)-a||_X= ||E(f|F_{\infty})-a||_X\ \forall a\in X\}=1$$

Since $X$ is separable and we can take $a=E(f(\o)|F_s)$ at every
$\o\in\O,$ we find that $E(f|F_s)\rightarrow E(f|F_{\infty})$ a.e.
as $s\to\infty.$

\end{pf}

\section{Martingale-ergodic and ergodic-martingale theorems}
In this section we prove norm as well as a.e. convergence for
vector valued martingale-ergodic and ergodic-martingale averages
with continuous parameter. In this section we consider only
regular martingales.

Following Kachurovskii \cite{Kach2}, we define martingale-ergodic
and ergodic-martingale averages as follows.

A \textit{martingale-ergodic} average is an average of the form
$\{E(A_tf|F_s)\}_{t>0, s\ge 0},$ where $E(\cdot|F_s)$ is the
conditional expectation operator and $A_tf$ is the ergodic average
while an \textit{ergodic-martingale} average is an average
$A_tE(f|F_s).$

Let us introduce the following notations
$$f_{\infty}(\o)=\lim\limits_{t\to\infty}A_tf(\o),$$
$$f^*(\o)=\lim\limits_{s\to\infty}E(f_{\infty}|F_s), \ \ f_*(\o)=\lim\limits_{t\to\infty}A_tE(f|F_{\infty})(\o).$$
The existence of above limit will be discussed below.

\begin{thm}
For $f\in L_p(\O, X),\ p\ge 1$ the following assertions hold.
\begin{enumerate}
    \item $$E(A_tf|F_s)\rightarrow f^*$$ in norm as $t,s\rightarrow \infty;$
    \item $$A_tE(f|F_s)\rightarrow f_*$$ in norm as $t,s\rightarrow \infty.$
\end{enumerate}

\end{thm}

\begin{pf} The idea is the same with real valued cases
\cite{Kach2}, \cite{pod1}.

Note that
$$||E(A_tf|F_s)-f^*||_p\le ||E(A_tf|F_s)-E(f_{\infty}|F_s)||_p+||E(f_{\infty}|F_s)-f^*||_p$$

The expression $||E(f_{\infty}|F_s)-f^*||_p$ converges due to
vector valued norm convergence theorem for continuous parameter
martingales \cite{Vaxan}.

Note that

$$||E(A_tf|F_s)-E(f_{\infty}|F_s)||_p=||E(A_tf-f_{\infty})|F_s)||_p\le ||A_tf-f_{\infty}||_p.$$

Since $||A_tf-f_{\infty}||_p$ convergent according to vector
valued ergodic theorem 2.1.5 \cite{kren}, then we get the
assertion (1).

Now, we prove the second part. According to Riesz convexity
theorem \cite{kren}, \cite{Phil} an $L_1-L_{\infty}$ contraction
is a contraction in $L_p$ norm. Therefore, we have the following
estimate

$$||A_tE(f|F_s)-f_*||_p\le ||A_tE(f|F_s)-A_tE(f|F_{\infty})||_p+||A_tE(f|F_{\infty})-f_*||_p\le$$
$$\le ||E(f|F_s)-E(f|F_{\infty})||_p+||A_tE(f|F_{\infty})-f_*||_p. $$

The norm $||E(f|F_s)-E(f|F_{\infty})||_p$ converges due to vector
valued norm convergence theorem for continuous parameter
martingales \cite{Vaxan}, and the norm
$||A_tE(f|F_{\infty})-f_*||_p$ from theorem 2.1.5 of \cite{kren}.

\end{pf}

We say that a linear operator $T$ in $L_1(\O, X)$ is
\textit{positively dominated} if there exists a positive linear
contraction $T'$ in $L_1,$ called a \textit{positive dominant} of
$T,$ such that
$$||Tf||_X\le T'(||f||_X).$$

Let us now provide some useful examples that we will use (see
\cite{SuchF}).

1. If $X=R$, then it is positively dominated by some positive
linear contraction on $L_1.$ For the vector valued $T$, a positive
dominant may not exist in general.

2. Let $\tau$ be a measure preserving transformation on $(\O
,\beta, \mu).$ Then the linear operator $T:L_1(\O, X)\rightarrow
L_1(\O, X)$ given by $Tf=f\circ\tau$ is said to be generated by
$\tau.$ $T$ is positively dominated by $T'$ with
$T'(||f||_X)=||f||_X\circ\tau.$

3. Assume that the Banach space $X$ has the Radon-Nikodym
property(the Banach space is said to have the Radon-Nykodim
property with respect to $(\O, \beta, \mu)$ if  any vector measure
$\phi:\beta\to X$ with finite variation, which is absolutely
continuous with respect to $\mu$ is just the integral of countable
valued function $f:\O\to X$ ). If $X$ is reflexive, then it has
the Radon-Nykodim property \cite{Vaxan}. Consider the conditional
expectation $E(f|F)$ with respect to $\sigma-$ subalgebra $F$ of
$\beta.$ For $f\in L_1(\O, X),$ the conditional expectation
$E(f|F)$ is Radon-Nikodym density with respect to the finite
measure $\mu$ on $F.$ Since $||E(f|F)||_X\le E'(||f||_X|F)$ a.e.
for all $f\in L_1(\O, X),$ where $E'(\cdot|F)$ is a conditional
expectation on $L_1,$ then the operator $E(\cdot|F)$ is positively
dominated by $E'(\cdot|F).$

We say that the flow $\{T_t, t\ge 0\}$ in $L_1(\O, X)$ is
positively dominated by the flow $\{P_t, t\ge 0\}$ in $L_1$ if for
any $f\in L_1(\O, X)$ and $t\ge0$ one has $||T_tf||_X\le
P_t(||f||_X)$ a.e. Now, we provide a.e. convergence theorem.

\begin{thm}  Let $X$ be a separable Banach space. Assume that $\{T_t, t\ge 0\}$ is
positively dominated by some semigroup $\{P_t, t\ge 0\}$ of
strongly continuous linear $L_1-L_{\infty}$ contractions. Then for
the function $f\in L_1(\O, X)$  the following assertions hold
true.

\begin{enumerate}
    \item If $\sup\limits_{t>0}||A_tf||_X\in L_1 $ (this holds for example, if $f\in L(\O, X)logL(\O, X)$ )  then for any $t>0, s\ge 0,$
    $E(A_tf|F_s)\rightarrow f^*$ a.e. as $t,s\rightarrow \infty.$
    \item If $\sup\limits_{s\ge 0}||E(f|F_s)||_X\in L_1,$ then $A_tE(f|F_s)\rightarrow
    f_*$ a.e. as $t,s\rightarrow \infty.$
\end{enumerate}

\end{thm}

\begin{pf} We prove the first assertion. Note that

$$||E(A_tf|F_s)-f^*||_X\le ||E(A_tf|F_s)-E(f_{\infty}|F_s)||_X+||E(f_{\infty}|F_s)-f^*||_X.$$

According to martingale convergence Theorem 2.4, the norm
$||E(f_{\infty}|F_s)-f^*||_X$ converges to $0$ a.e. as
$s\to\infty.$

(Let $0<t_1\le t.$) Further, since conditional expectation
operator is positively dominated, then

$$||E(A_tf|F_s)-E(f_{\infty}|F_s)||_X=||E(A_tf-f_{\infty}|F_s)||_X\le E'(||A_tf-f_{\infty}||_X|F_s)\le E'(h_{t_1}|F_s),$$
where $h_{t_1}(\o)=\sup\limits_{t\ge
t_1}||A_tf(\o)-f_{\infty}(\o)||_X$ and $E'$ is a positive dominant
of $E.$ Due to the condition of the theorem, we have  $h_{t_1}\in
L_1$ and $h_{t_1}\to 0$ a.e. from  Theorem 2.1. Now applying first
part of Theorem 1.2, $E'(h_{t_1}|F_s)\rightarrow 0,$ a.e. as
$t_1\to\infty.$ Therefore, we have
$||E(A_tf|F_s)-E(f_{\infty}|F_s)||_X\rightarrow 0,$ a.e. Hence,
$||E(A_tf|F_s)-f^*||_X\rightarrow 0$ a.e. as $t,s\to\infty.$

Now we prove the second part. We have

$$||A_tE(f|F_s)-f_*||_X\le ||A_tE(f|F_s)-A_tE(f|F_{\infty})||_X+||A_tE(f|F_{\infty})-f_*||_X$$

The norm $||A_tE(f|F_{\infty})-f_*||_X$ is a.e. convergent due to
Theorem 2.1.

We have the following
$$||A_tE(f|F_s)-A_tE(f|F_{\infty})||_X=||\frac1t\int\limits_0^tT_{\tau}(E(f|F_s)-E(f|F_{\infty}))d\tau||_X\le$$
$$\le \frac1t\int\limits_0^t||T_{\tau}(E(f|F_s)-E(f|F_{\infty}))||_Xd\tau\le $$$$\le \frac1t\int\limits_0^tP_{\tau}\big(||(E(f|F_s)-E(f|F_{\infty}))||_X\big)d\tau= $$
$$=A'_t\big(||(E(f|F_s)-E(f|F_{\infty}))||_X\big)$$

where $P_{t}$ is a positive dominant of $T_t$ for each $t$ and
$A'_tf=\frac1t\int\limits_0^tP_{\tau}fd\tau.$ According to our
assumption, the flow $\{P_t, t\ge 0\}$ is strongly continuous
semigroup.

Note that the real valued function
$h_s(\o)=||(E(f(\o)|F_s)-E(f_{\infty}(\o)|F_s))||_X$ is integrable
according to the conditions of theorem. Moreover, according to the
martingale convergence Theorem 2.3 $h_s(\o)\rightarrow 0$ a.e. as
$s\to\infty.$

Now applying second part of Theorem 1.2, we get
$A'_t(h_s)\rightarrow 0$ a.e. as $t,s\to\infty.$ Therefore,
$||A_tE(f|F_s)-A_tE(f|F_{\infty})||_X\rightarrow 0$ a.e. as
$s,t\to\infty.$

\end{pf}

\textbf{Remark}. When we consider real valued functions, that is
when $X=R,$ then for any semigroup $\{T_t, t\ge 0\}$ of linear
$L_1-L_{\infty}$ contractions there always exists a semigroup
$\{P_t, t\ge 0\}$ of positive linear $L_1-L_{\infty}$ contractions
such that $|T_tf|\le P_t|f|$ a.e. However, in vector valued
positive dominant semigroup may not exist in general. It is also
known that $\{T_t, t\ge 0\}$ is not positively dominated by its
linear modulus \cite{Sato3}. Therefore in the above theorem,
despite real valued case, we need an additional assumption that
$\{T_t, t\ge 0\}$ should be positively dominated by  $\{P_t, t\ge
0\}.$  Of course one can ask to provide the above theorems without
this condition, but we fail to answer to this question.

The following theorem is  dominant and maximal inequalities for
martingale-ergodic processes.

\begin{thm} Under the assumption of Theorem 3.2, $f\in L_p(\O, X),\ p>1,$  $\sup\limits_{t>0}||A_tf||_X\in L_1$ and  $F_s\downarrow F,\ s\to\infty$
then the following assertions hold true.

\begin{enumerate}
    \item
    $$||\sup\limits_{t,s}||E(A_tf|F_s)||_X||_p\le \big(\frac p{p-1}\big)^2||f||_p,$$

    \item $$\mu\big\{\sup\limits_{t,s}||E(A_tf|F_s)||_X\ge \varepsilon\big\}\le \frac p{p-1}\frac{||f||_p}{\varepsilon}$$
\end{enumerate}

\end{thm}

\begin{pf} We first prove the dominant inequality. Note that the
conditional expectation operator is positively dominated, then

$$||\sup\limits_{t,s}||E(A_tf|F_s)||_X||_p\le||\sup\limits_{t,s}E'(||A_tf||_X|F_s)||_p\,$$

where $E'$ is a positive dominant of $E.$

Since $\{T_t, \ t\ge 0\}$ is positively dominated by $\{P_t, \
t\ge 0\}$, then

$$||A_tf||_X=||\frac1t\int\limits_0^tT_tfd\tau||_X\le \frac1t\int\limits_0^t||T_tf||_Xd\tau\le$$
$$\le \frac1t\int\limits_0^tP_t(||f||_X)d\tau=A'_t(||f||_X).$$

Since $E$ is positively dominated by $E'$ and $A_t$ by $A'_t,$
then we have the following inequality.
$$\sup\limits_{t,s}E'(||A_tf||_X|F_s)||_p\le \sup\limits_{t,s}E'(A'_t(||f||_X)|F_s)||_p,$$

 Since the flow
$P_t$ is a strongly continuous semigroup, applying Theorem 3 of
\cite{pod1} for the process $E'(A'_t(||f||_X)|F_s)$, we get

$$||\sup\limits_{t,s}E'(A'_t(||f||_X)|F_s)||_p\le \big(\frac p{p-1}\big)^2||f||_p.$$

The above chain of inequalities imply part (1) of the theorem.

Now we prove part (2). Since the operator $E$ is positively
dominated by some $E'$ , then we have the following inequalities

$$\mu\big\{\sup\limits_{t,s}||E(A_tf|F_s)||_X\ge \varepsilon\big\}\le \mu\big\{\sup\limits_{t,s}E'(||A_tf||_X|F_s)\ge \varepsilon\big\}\le$$
$$\le \mu\big\{\sup\limits_{t,s}E'(A'_t||f||_X|F_s)\ge \varepsilon\big\}$$
where
$A'_t(||f||_X)=\frac1t\int\limits_0^tP_{\tau}(||f||_X)d\tau.$ Now
applying second part of Theorem 3 of \cite{pod1}, for the process
$E'(A'_t(||f||_X)|F_s)$, we get

$$ \mu\big\{\sup\limits_{t,s}E'(A'_t||f||_X|F_s)\ge \varepsilon\big\}\le \frac p{p-1}\frac{||f||_p}{\varepsilon}.$$

Hence (2) is proved.

\end{pf}

Now, we provide dominant and maximal inequalities for
ergodic-martingale average.

\begin{thm} Under the assumption of Theorem 3.2, $f\in L_p(\O, X),\ p>1,$ $\sup\limits_{s\ge0}||E(f|F_s)||_X\in L_1$ and  $F_s\downarrow F,\ s\to\infty.$
and the following assertions hold true.

\begin{enumerate}
    \item
    $$||\sup\limits_{t,s}||A_tE(f|F_s)||_X||_p\le \big(\frac p{p-1}\big)^2||f||_p,$$

    \item $$\mu\big\{\sup\limits_{t,s}||A_tE(f|F_s)||_X\ge \varepsilon\big\}\le \frac p{p-1}\frac{||f||_p}{\varepsilon}.$$
\end{enumerate}

\end{thm}

 This theorem can easily be proven using Theorem 4 of
\cite{pod1} and the way of proof of Theorem 3.3. So we omit
details.

It is known that in $L_1,$ the condition of integrability of
supremum can not be omitted in all unified theorem \cite{ArgRos}.
The following theorem is given without this assumption, but the
conditional expectation operator and ergodic average should
commute.

\begin{thm} Let $F_s\downarrow F,\ s\to\infty$ and  $T_t$ be a semigroup of strongly continuous
measure preserving transformation and $T_tE(f|F_s)=E(T_tf|F_s),$
for all $t,s\ge 0.$ Then for any $f\in L_1(\O, X),$ the averages
$A_tE(f|F_s)$ and $E(A_tf|F_s)$ converge a.e. as $t,s\to\infty.$

\end{thm}

\begin{pf} The idea is almost the same as  Theorem 4 of
\cite{pod2}.

Let $n=[t],$ then $n=t+\a,$ where $0\le\a<1.$ For any $t>0,\ s\ge
0$ we have

$$A_tE(f|F_s)=\frac1t\int\limits_0^tT_{\tau}E(f|F_s)d\tau=\frac1t\int\limits_0^nT_{\tau}E(f|F_s)d\tau+\frac1t\int\limits_n^{n+\a}T_{\tau}E(f|F_s)d\tau=$$
$$=\frac1t\sum\limits_{k=0}^{n-1}\int\limits_k^{k+1}T_{\tau}E(f|F_s)d\tau+\frac1t\int\limits_n^{n+\a}T_{\tau}E(f|F_s)d\tau=$$
$$=\frac1t\sum\limits_{k=0}^{n-1}\int\limits_0^1T_{\tau+k}E(f|F_s)d\tau+\frac1t\int\limits_0^{\a}T_{\tau+n}E(f|F_s)d\tau=$$
$$=\frac1t\sum\limits_{k=0}^{n-1}T_k\int\limits_0^1T_{\tau}E(f|F_s)d\tau+\frac1tT_n\int\limits_0^{\a}T_{\tau}E(f|F_s)d\tau=$$
$$=\frac1t\sum\limits_{k=0}^{n-1}(T_1)^kE(A_1f|F_s)d\tau+\frac{\a}t(T_1)^nA_{\a}E(f|F_s)=$$
$$=\frac nt[S_n(T_1)E(g_1|F_s)+\frac{\a}n(T_1)^nA_{\a}E(f|F_s)],$$
where $g_1=A_1f$ and $S_n(T)f=\frac
1n\sum\limits_{i=0}^{n-1}T_if.$

Now let us estimate the expressions $S_n(T_1)E(g_1|F_s)$ and
$\frac{\a}n(T_1)^nA_{\a}E(f|F_s).$ Evidently, the former is a.e.
convergent. If $P^1$ and $E'$ be positive dominants of $T_1$ and
$E$ respectively, then the latter converges a.e. since
$$||\frac{\a}n(T_1)^nA_{\a}E(f|F_s)||_X\le \frac 1n (P^1)^nE'(A_1||f||_X|F_s)=$$
$$=\frac{n+1}nS_{n+1}(P^1)E'(||f||_X|F_s)-S_n(P^1)E'(||f||_X|F_s)\rightarrow 0$$
a.e. as $s,n\to\infty$ from Theorem 1.2.

\end{pf}

\end{document}